
\font\eightrm=cmr8         
\font\eightit=cmti8        
\font\eightmi=cmmi8        

\def\cite#1{[#1]}             
\def\frac#1#2{{#1\over #2}}   

\def\N{{\bf N}}               
\def\R{{\bf R}}               

\def\<#1>{\langle#1\rangle}   
\def\TDS{{\eightrm TDS}}


\outer\def\section#1. #2\par{
      \bigskip\bigskip       
      \message{#1. #2}
      \leftline{\bf#1. #2}
      \nobreak\smallskip     
      \noindent}             

\outer\def\varsection#1\par{
      \bigskip\bigskip       
      \message{#1 }
      \leftline{\bf#1 }
      \nobreak\smallskip     
      \noindent}             

\def\declare#1. #2\par{
      \medskip\noindent     
      {\bf#1.}\rm           
      \enspace\ignorespaces 
      #2\par\smallskip}     



\def\cite#1{{\rm[#1]}}        
\def\refno#1. #2\par{\smallskip\item{[#1]} #2\par}

\def\DonatoGeometrie{1}
\def\DonatoVirasoro{2}
\def\SouriauAlgoritme{3}
\def\Sou{4}
\def\Torre{5}
\def\TorreTangen{6}
\def\TorreBan{7}


\vsize= 51truepc            
\hsize= 36truepc            
\voffset= 14truept          

\hoffset= 18truept          

\newtoks\rightheadtext      
\newtoks\leftheadtext       

\headline={\ifnum\pageno>1  
            \ifodd\pageno   
             \hfil \eightrm \the\rightheadtext \hfil
             \llap{\tenrm\folio}%
            \else \rlap{\tenrm\folio}
             \hfil \eightrm \the\leftheadtext \hfil \fi
           \else\hfil \fi}
\footline={\hfil}           

\rightheadtext{THE DEFINITION OF A DIFFERENTIAL FORM}

\leftheadtext{CARLOS A. TORRE}


\topglue 1.5cm

\centerline{\bf  THE DEFINITION OF A DIFFERENTIAL FORM}
\centerline{\bf ON A TANGENT STRUCTURE}

\bigskip

\centerline{\sl Carlos A. Torre}
\bigskip
\centerline{\sl Escuela de Matem\'atica,
                 Universidad de Costa Rica,}
\centerline{\sl San Jos\'e, Costa Rica}
\bigskip

{\narrower
\eightrm
\baselineskip=9.5pt
\textfont1=\eightmi
\noindent
{\eightit ABSTRACT}.
Three definitions of a differential form on a tangent structure
are considered. It is proved that the (covariant) definition given by Souriau
(as a collection of forms indexed by the plaques) is equivalent
to a smooth section of the corresponding vector bundle if the space
does not have transverse points.

\par
\smallskip}

\bigskip

\section 1. Introduction.

A diffeology on any set $X$ was defined by Souriau [3,4] as a collection
of plaques $p\colon U\subset \R^n \to X$, for any $n\in \N$ which covers
$X$ and is closed under composition with smooth maps 
$\phi \colon V\subset \R^m \to \R^n$ where $U, V$ are open sets. In this
context
differential forms were defined as a collection of forms on open
subsets of $\R^n$, indexed by the plaques. They were used to construct
a picture of quantization on a class of coadjoint orbits of 
diffeomorphisms groups [1,2]. A tangent structure was defined in [7]
for coadjoint orbits of diffeomorphisms groups, allowing a geometric
definition of a differential form as a section of the corresponding
vector bundle [5,6], and in this way allowing also the construction of
a Poisson structure and precuantization on coadjoint orbits 
of diffeomorphisms groups. In this paper we consider the relation between
the two definitions of a differential form mentioned above, and the
more algebraic one as a multilinear alternating map on the 
$C^{\infty}(X)$ module of vector fields on $X$.

In section 2 we review four definitions of a differential form on
a finite dimensional manifold and recall the proof about 
the equivalence between them.
In section 3 we recall the definition of a tangent structure and
three definitions of a differential form on it: the global definition
given algebraically, the local definition, as a collection of
forms indexed by plaques, as defined by Souriau, 
and as a smooth section on the vector bundle
of exterior forms. It is proved that the last two are equivalent if
the space does not have transverse points.

\section 2. A review of differential forms on manifolds.

An {\sl exterior k-form} on a vector space $V$ is an alternating
k-linear map 
$$\omega \colon V^k \to \R$$ 
The set of them is 
denoted $\Lambda^kV^*$, forming a vector space with the ordinary
scalar multiplication and sum of functions and the set
$\Lambda^{\cdot}V^* := \cup_{k\in \N}\Lambda^kV^*$ is an associative
algebra with the {\sl wedge product}:
$$\omega^k \wedge \omega^l(\xi_1, \cdots , \xi_{k+l}) :=
\sum (-1)^{\sigma}\omega^k(\xi_{i_1}, \cdots ,\xi_{i_k})
\omega^l(\xi_{j_1}, \cdots ,\xi_{j_l})$$
where $\omega^i\in \Lambda^iV^*$ for $i=k,l$ and the sum is over all
permutations $(i_1, \cdots ,i_k, j_1, \cdots ,j_l)$ of
$(1, \cdots ,k+l)$ and $\sigma$ is the order of the permutation. This
operation is skew commutative:
$$\omega^k \wedge \omega^l = (-1)^{kl}\omega^l \wedge \omega^k$$
In particular if $k=l=1$ then
$$\omega_1 \wedge \omega_2(\xi, \eta) = \omega_1(\xi)\omega_2(\eta) -
\omega_2(\xi)\omega_1(\eta)$$
and more generally
$$\omega_1 \wedge \cdots \wedge \omega_k(\xi_1, \cdots ,\xi_k) =
det(\omega_j(\xi_i))$$
In $V= \R^n$ every exterior 1-form $\omega$ is a scalar multiple of a 
projection along some vector $y\in V$: $\omega(\xi) = (\xi, y)$. In 
particular the projection along the i-th coordinate vector is denoted
$x_i$. More generally, if $P\colon V \to F$ is a projection along a
k-dimensional subspace $F$ with ortonormal base
$\{ \beta_1, \cdots \beta_k \}$, then the map
$$\omega^k_P(\xi_1, \cdots ,\xi_k) = det(a_{ij})$$
where $P(\xi_i)= \sum_{j=1}^ka_{ij}\beta_j$ defines a k-form on $V$. 
This is the oriented k-volume of the projection along $F$ of the 
parallelepiped generated by $\xi_1, \cdots \xi_k$.
Conversly, using linearity it is found that every exterior k-form is a
linear combination of forms $\omega^k_P$ where $P$ runs over all 
projections along the k-dimensional coordinate subspaces and therefore 
the set
$$\{ x_{i_1} \wedge \cdots \wedge x_{i_k} / \{ i_1, \cdots i_k\}
\subset \{ 1, \cdots n\} \}$$
is a base of $\Lambda^kV^*$.
In particular if $k=n$, every exterior n-form is a scalar
multiple of 
$$\nu(\xi_1, \cdots \xi_n) = det(\xi_{ij})$$
where $\xi_i = \sum_{i=1}^n\xi_{ij}e_j$ and 
$\{ e_1, \cdots e_n\}$ is the standard base. This is the oriented
volume of the parallelepiped with edges $\xi_1, \cdots \xi_n$.
Every k-form with $k>n$ is zero.

If $\omega_1, \omega_2$ are projections along $v_1, v_2$ then
$\omega_1 \wedge \omega_2(\xi_1, \xi_2)$ is the oriented area of the
projection of the parallelepiped generated by $\xi_1, \xi_2$ over the
plane generated by $v_1, v_2$.

The {\sl pullback} of $\omega \in \Lambda^kV^*$ by a linear map
$L\colon W \to V$ is the exterior k-form $L^*\omega \in \Lambda^kW^*$
defined by
$$L^*\omega (\xi_1, \cdots ,\xi_k) = \omega (L\xi_1, \cdots L\xi_k)$$
The pullback defines an algebra homomorphism and preserves composition
of linear maps (in the appropiate order).

A {\sl differential k-form} on an open set $U\subset \R^n$ is an 
alternating $C^{\infty}(U)$ k-linear map
$$\omega \colon [\Gamma(U)]^k \to C^{\infty}(U)$$
where $\Gamma(U)$ denotes the $C^{\infty}(U)$ -module of smooth
vector fields on U. Therefore, each $x\in U$ defines 
$\omega_x \in \Lambda^kV^*$ by
$$\omega_x(\xi_1, \cdots ,\xi_k) := \omega(\bar{\xi_1}, \cdots ,\bar{\xi_k})
(x)$$
where $\bar{\xi_1}, \cdots \bar{\xi_k}$ are smooth vector fields such
that $\bar{\xi_i}(x)= \xi_i$ (they exist:
$\bar{\xi_i}(x) = \sum \xi_i{\partial \over \partial x_i}$).

$\omega_x$ is well defined: if $\xi(x_0)=0$, then using the
$C^{\infty}(U)$-linearity and the expression
$\xi = \sum_jb_j{\partial \over \partial x_j}$, with 
$b_j\in C^{\infty}(U)$ and $b_j(x_0)=0$ it follows that
$\omega(\xi)(x_0)=0$.

Conversely, by evaluation on the vector fields 
${\partial \over \partial x_j}$, it follows that an exterior k-form
$$\omega_x = \sum_{\vert I\vert =k} 
a_{i_1\cdots i_k}(x)x_{i_1}\wedge \cdots x_{i_k}$$
for each $x\in U$ defines a differential k-form on $U$ iff each
$a_{i_1\cdots i_k}\in C^{\infty}(U)$
In particular, the k-forms defined by $x_{i_1}\wedge \cdots \wedge x_{i_k}$
for each $x\in U$ is denoted $dx_{i_1}\wedge \cdots \wedge dx_{i_k}$, since
$dx_i$denoted the 1-form dual of the vector field 
${\partial \over \partial x_i}$ and $dx_i$ is also the differential of 
the projection along the i-axis (denoted $x_i$ previously).
Therefore, we have two equivalent definitions of a k-form on $U$. 

If $f\colon V\subset \R^m \to U$ is smooth, the {\sl pullback} of a k-form
$\omega$ on $U$ is defined using the differential of $f$:
$$f^*\omega(\xi_1, \cdots \xi_k)(x)=
\omega_{f(x)}(df(\xi_1)(x), \cdots ,df(\xi_k)(x))$$

\declare Definition 2.1.
A k-form on a smooth manifold $M$ given by an atlas
${\cal A} = \{ (\alpha_i, U_i) / i\in I \}$
is defined in four different ways (al least):
\item{1.}
It is an exterior k-form $\omega_x$ on $T_xM$ for each $x\in M$ such that
if $\xi_i\in \Gamma (M)$ for $i=1, \cdots k$, then
$$\omega(\xi_1, \cdots ,\xi_k)(x) := \omega_x(\xi_1(x), \cdots ,\xi_k(x))$$
is smooth.
\item{2.}
It is a $C^{\infty}(M)$ k-linear alternating map
$$\omega \colon [\Gamma(M)]^k \to C^{\infty}(M)$$
\item{3.}
It is a collection $(\omega_i)_{i\in I}$ of k-forms, where
$\omega_i$ is a k-form on $U_i\subset \R^n$ such that 
$$\omega_j = (\alpha_i^{-1}\circ \alpha_j)^*\omega_i$$
for each $i,j\in I$.
\item{4.}
It is a smooth section on the vector bundle of exterior k-forms on $M$.

The four definitions are equivalent: $1\Rightarrow 2$ follows from 
evaluation at a point. $2\Rightarrow 1$ follows as follows: from a vector
$\xi_i = \sum \xi_{ij}{\partial \over \partial x_j}\in T_{x_0}M$ construct a
vector field $\bar{\xi_i}(x) = \sum \xi_{ij}{\partial \over \partial x_j}$
in a neighborhood of $x_0$ and extend globally through multiplication by 
a smooth function $g$ such that $x_0\in V<g<U$ (that is: $Supp g \subset U$,
$g=0$ on $V$ and $V$ is a neighborhood of $x_0$ such that $\bar V\subset U$)
and defining $\bar{\xi_i}(x)=0$ outside $U$. 
 Then $\omega_x$
is defined using the vector fields $\bar{\xi_i}$. In order to prove that
the definiton is independent of this particular choice of $\bar{\xi_i}$,
assume that $\xi(x_0)=0$, take $U,g$ as before and a local expression
$\xi(x)=\sum b_j(x){\partial \over \partial x_j}$ on $U$. Then
$$\omega(\xi)(x_0) = g(x_0)\sum b_j(x_0){\partial \over \partial x_j}=0$$ 
$1\Rightarrow 3$ is obtained as follows: define $\omega_i= \alpha_i^*\omega$,
use $\alpha_j = \alpha_i\circ (\alpha_i^{-1}\circ \alpha_j)$ and the fact
that the pullback preserves composition and $d\alpha_i$ preserves smoothness
of vector fields.
In order to prove $3\Rightarrow 1$, for each $x\in M$ take $i\in I$ such 
$x\in U_i$ and define $\omega_x = (\alpha_i^{-1})^*\omega_i(x)$ which
does not depend on the chart chosen.
Finally $1\iff 4$: the vector bundle is given by 
$\wedge^kT^*M \to M$ where
$$\wedge^kT^*M = \{ (x,\xi) / x\in M, \xi \in \wedge^kT^*_xM \}$$
$\pi(x,\xi)= x$ and with the atlas
$$\{ (\pi^{-1}(U_1), (\alpha_i\times id)\circ \phi_i) / i\in I \}$$
where $\phi_i\colon \pi^{-1}(U_i) \to U_i\times \R^k$
is given by
$\phi_i(x,\xi) := (x,\xi_{i_1\cdots i_k}(x))$ such that
$\xi = \sum \xi_{i_1\cdots i_k}dx_{i_1}\wedge \cdots dx_{i_k}$.
Now, given $\omega$ under definition 1, define $\bar{\omega}$ under
definition 4 (and viseversa) by
$$\bar{\omega}(x)(\xi_1(x), \cdots \xi_k(x) =
(x, \omega_x(\xi_1(x), \cdots ,\xi_k(x))$$
The smoothness of $\bar{\omega}$ and $\omega$ are equivalent using the
definition of smoothness with the local expresssion of these forms.
This proves the equivalence of the four definitions.

The space of k-forms on $M$ is denoted $\epsilon^k(M)$. It is a
$C^{\infty}(M)$-module and $\epsilon^{\cdot}(M) :=\cup_{k\geq 0}\epsilon^k(M)$
is an associative algebra with the wedge product, where
$\epsilon^0(M)$ is defined as $C^{\infty}(M)$.

For example, if $f\colon M\to \R$ is a smooth function, then $df$ defined
by $df_x(\xi_x) = {d\over dt}\vert_0d(\alpha(t))$ where 
$\alpha^{\prime}(0)=\xi_x$ defines a 1-form. For every 1-form $\omega$ on a 
Riemannian manifold $M$ there exists a smooth vector field $\eta$ such
that $\omega(\xi)(x) = (\xi_x, \eta_x)$. In particular, if $\omega = df$
then $\eta$ is called the gradient of $f$.

\section 3. Tangent structures.

An {\sl n-plaque} on a set $X$ is a map $p\colon U\subset \R^n \to X$
with $U$ open. A {\sl diffeology on X} is a collection $P(X)$ of n-plaques
for each $n\in \N$ which covers $X$ and such that if $p\in P(X)$ and
$\phi \colon U\subset \R^m \to \R^n$ is smooth, then $p\circ \phi\in P(X)$.
A {\sl tangent structure} on $(X, P(X))$ is a collection of equivalence
relations $\sim_F^n$ on the set of n-plaques at $F\in X$ (that is,
$p\in P(X), p(0)= F$) denoted $P^n_F(X)$, for each $n\in \N, F\in X$,
satisfying the following two consistency conditions:
\item{1.}
If $p_1 \sim_F^n p_2$ and $\phi \colon U\subset \R^m \to \R$, then
$p_1\circ \phi \sim_F^m p_2\circ \phi$.
\item{2.}
If $p\colon U\to X$ is a plaque and $V\subset U$ is an open neighborhood
of $0$, then $p\sim p\vert_V$. 

The class of $p$ is denoted $[p]$ or
$[p(t)]_t$. A diffeological space with a tangent structure is called
a $\TDS$.
The tangent structure is called {\sl linear} if $P_F^n(X)/\sim_F^n$
carries a vector space structure for each $F\in X, n\in \N$ (this set is
called the {\sl n-th tangent space at $F$} and is denoted $T_F^nX$)
satisfying
\item{a.}
If $p_{12}\in [p_1]+[p_2]$ and $\phi \colon U\subset \R^m \to \R^n$ is
smooth, then 
$$p_{12}\circ \phi \in [p_1\circ \phi] + [p_2 \circ \phi]$$
\item{b.}
If $p\in c[p_1] \Rightarrow p\circ \phi \in c[p_1\circ \phi]$.

The tangent structure is called {\sl continuous} if given two (n+m)-plaques
$p_i(r,s), i=1,2$ such that $p_1(r,0)= p_2(r,0)$ for each $r$, there
exists a (n+m)-plaque $p_{12}(r,s)$ such that
$$[p_{12}(r,s)]_s = [p_1(r,s)]_s + [p_2(r,s)]_s \ \ \ \forall r$$
If $(X, P(X), \sim), (Y, P(Y), \equiv)$ are two $\TDS$, a map
$f\colon X \to Y$ is called {\sl smooth} if it preserves plaques, directions,
and the vector space structure, that is:
\item{1.}
$p\in P(X) \Rightarrow f\circ p\in P(Y)$
\item{2.}
$p_1 \sim_F^n p_2 \Rightarrow f\circ p_1 \equiv_{f(F)}f\circ p_2$
\item{3.}
$d^kf\colon T^kX \to T^kY$ given by $[p] \to [f\circ p]$ is linear
if the tangent structure is linear.

The set of smooth maps is denoted $C^{\infty}(X,Y)$ and the set of 
maps satisfying only $1.$ is denoted $c^{\infty}(X,Y)$.

Any manifold $M$ has a {\sl standard $\TDS$} for which $P(M)$ is formed by
the smooth maps $p\colon U\subset \R^n \to M$ (which contains the atlas)
and 
$$p_1 \sim_F^n p_2 \iff D^ip_1\vert_0 = D^ip_2\vert_0 \ \ \forall
i=0, 1, \cdots ,n.$$
The tangent spaces and the smooth maps with this $\TDS$
coincide with those given by
the manifold. In particular the standard manifold structure on $\R^n$ 
provides a standard $\TDS$.
Using this, any subalgebra $\cal F$ of $c^{\infty}(X)$ on any diffeology
$(X, P(X))$ allows the construction of a tangent structure on
$(X, P(X))$, called the {\sl standard $\TDS$ defined by $\cal F$}:
$$p_1\sim_F^n p_2 \iff D^if\circ p_1\vert_0 = D^if\circ p_2\vert_0 \ \ 
\forall f\in {\cal F}$$
in which we are using the standard $\TDS$ on $\R$.

The tangent spaces defines the {\sl tangent bundles} $T^nX$. A 
{\sl vector field} is a section of this bundle such that 
$\forall f\in C^{\infty}(X, \R)$ the map
$df(\xi)\colon X \to \R$ given by $F\to df(\xi(F))$ is smooth.
The set of smooth vector fields is denoted $\Gamma(X)$. In every linear
$\TDS$ the set $\Gamma(X)$ is a $C^{\infty}(X)$ module. 

A vector field
is called {\sl locally integrable} if for every n-plaque
$p\colon U\subset \R^n \to X$ and every $r_0\in U$ there exists a 
neighborhood $V$ of $r_0$ and a (n+1)-plaque $q$ such that 
$q(r,0)= p(r)$ and
$$\xi(p(r))= [q(r,t)]_t \ \ \forall r\in V_0$$
Notice that every smooth vector field on a manifold is locally 
integrable. If the tangent structure is linear and continuous, the set
of locally integrable vector fields is a $C^{\infty}(X)$ module,
denoted $\Gamma_I(X)$.

A point $F\in X$ is called {\sl weakly (strongly resp.) transverse
of degree $(n,m)$} if there exists $p_1\colon U_1 \subset \R^n \to X$,
$p_2\colon U_2\subset \R^m \to X$ at $F$ such that there does not exist
$ V_1 \subset U_1, V_2\subset U_2$ open at $F$ and a (n+m)-plaque
$q\colon V_1\times V_2 \to X$ satisfying $q(r,0)=p_1(r), \forall r$
and $q(0,s)=p_2(s), \forall s$
(or $[q(r,0)]_r= [p_1(r)]_r, [q(0,s)]_s$ respectively)

\declare Examples 3.1.

\item{[1]}
Manifolds modelled on a locally convex vector space $E$ do not have transverse
points: if $\{ (\alpha_i, U_i) / i\in I \}$ is an atlas of $M$ and
$p_1, p_2$ are n, m-plaques respectively, at $F\in M$ then
$p_i = \alpha_i\circ \phi_i, i=1,2$ where $\phi_i$ is smooth. We may
assume that 
$\alpha_1= \alpha_2$, $\phi_1\colon V_1\subset \R^n \to U\subset E$,
$\phi_2\colon V_2\subset \R^m \to U$. Consider $V_i^{\prime}\subset V_i$
such that 
$$\phi_1(V_1^{\prime}) + \phi_2(V_2^{\prime})\subset U$$
Let $\phi \colon V_1^{\prime} \times V_2^{\prime} \to U$ defined by
$\phi(u,v)= \phi_1(u)+\phi_2(v)$, then $q= \alpha_1\circ \phi$
satisfies the required condition.

\item{[2]}
Let $X$ be the group of diffeomorphisms of a finite dimensional manifold
$M$. Let $P(X)$ be the set of plaques $p\colon U\subset \R^n \to X$ such 
that $\bar p\colon U\times M \to M$ given by $\bar p(r,m)=p(r)(m)$ is 
smooth and consider the standard tangent structure defined by
$c^{\infty}(X)$. Then $X$ does not have transverse points: take
$q(r,s)= F\circ p_1(r)p_2(s)\circ F^{-1}$.

\item{[3]}
Let $X= \R^2$, let $P(X)$ be the set of smooth maps 
$p\colon U\subset \R^n \to X$ whose image is a line (under the standard
$\TDS$), then $c^{\infty}(X)$ consists of maps which are smooth along 
lines. With the standard $\TDS$ the tangent structure is linear and
continuous, each $T_FX$ is $\R^2$. The only locally integrable
vector field is the zero vector field. Notice that every point is
strongly transverse of degree $(1,1)$. More generally we consider a 
manifold $M$ (or a collection of manifolds $(M_i)_{i\in I}$)
and a set of smooth maps $p\colon U\subset \R^n \to M$, for different
$n\in \N$, then consider the diffeology generated by these maps. For 
example, we take $M= S^2$ and consider the maps whose image is a segment
of a parallel (notice that the intersection of 2 plaques is either a point
of another plaque). The poles are strongly transverse points of
degree $(1,1)$. A vector field is locally integrable $\iff$ it is zero
at poles. However, the union of two tangent surfaces 
(in this context) at only one
point $F$ do not have transverse points and every vector field is locally
integrable.

Notice that every strongly transverse point is weakly transverse. Every
smooth vector field is zero at a strongly transverse point.

\declare Definition 3.1.

Given a linear $\TDS$ $(X, P(X), \sim)$, a k-form on $X$ is defined in
one of the following ways:
\item{1.}
It is an exterior k-form $\omega_F, \forall F\in X$ such that if
$\xi_1, \cdots ,\xi_k$ are vector fields then
$$\omega(\xi_1, \cdots ,\xi_k)(F) := \omega_F(\xi_1(F), \cdots ,\xi_k(F))$$
is a smooth map.
\item{2.}
It is a $C^{\infty}(X, \R)$ alternating k-linear map
$$\omega \colon (\Gamma(X))^k \to C^{\infty}(X, \R)$$
\item{3.}
It is a collection $(\omega_p)_{p\in P(X)}$ where $\omega_p$ is a k-form
on $U\subset \R^n$ if $p\colon U\to X$ such that
\itemitem{a.}
$\omega_{p\circ \phi} = \phi^*\omega_p, \forall \phi \colon U^{\prime}
\subset \R^m \to \R^n$ smooth.
\itemitem{b.}
If $p_1, p_2$ are two plaques satisfying $p_1(r_1)= p_2(r_2)=F$ which
are tangent at $F$ along directions $v_1, \cdots ,v_k\in \R^n$ (this means
that
$[p_1(tv_i+r_1)]_t = [p_2(tv_i+r_2)]_t$ for all $i=1, \cdots ,k$)
then
$$\omega_{p_1(r_1)}(v_1, \cdots ,v_k)=
\omega_{p_2(r_2)}(v_1, \cdots ,v_k)$$
and
$$d\omega_{p_1(r_1)}(v_1, \cdots ,v_{k+1})=
d\omega_{p_2(r_2)}(v_1, \cdots ,v_{k+1})$$
for all $v_{k+1}\in \R^n$. This is called the {\sl tangent condition}
for further reference.

The set of k-forms under definition $i$ is denoted $\epsilon^k_i(X)$
for $i=1,2,3$. Each of them is a $C^{\infty}(X)$ module and
$\epsilon^{\cdot}_i(X)$, is also an associative algebra
with the wedge product defined as usual under each of the three
definitions.

The pullback of $\omega \in \epsilon_1^k(X_2)$ under a smooth map
$h\colon X_1 \to X_2$ is defined by
$$(h^*\omega)_F(\eta_1(F), \cdots ,\eta_k(F))
=\omega_{\eta(F)}(dh(F)\eta_1(F), \cdots ,dh(F)\eta_k(F))$$
If $\omega \in \epsilon_3^k(X_2)$ then define
$$(h^*\omega) = (\omega_{p\circ h})_{p\in P(X_1)}$$
This collection satisfies the tangent condition.
If $\omega \in \epsilon_2^k(X_2)$ we need $dh(F)$ to be surjective 
$\forall F\in X$ and define
$$h^*\omega(\eta_1, \cdots ,\eta_k)= \omega(dh(\eta_1), \cdots,
dh(\eta_k))\circ h$$

\declare Proposition 3.1.

If $(X, P(X), \sim)$ does not have transverse points then
$\epsilon_1^k(X)= \epsilon_3^k(X)$ for any $k\in \N$

{\sl Proof:}
First we define a map $\Psi \colon \epsilon_1^k(X) \to \epsilon_3^k(X)$
by 
$$\Psi (\omega) = (p^*\omega)_{p\in P(X)}$$
This map is linear and $\Psi(f\omega)= f\Psi(\omega)$ if $f\in C^{\infty}(X)$
We prove that $\Psi$ is 1-1:
Assume that $\Psi(\omega)= 0$ and that $\omega_F \neq 0$ for some $F\in X$,
then $\exists v_1, \cdots ,v_k\in \underline{T_FX}$ such that 
$\omega_F(v_1, \cdots ,v_k)\neq 0$. Let $p_1, \cdots p_k$ be 1-plaques such
that $[p_i]= v_i$. Since $F$ is not strongly transverse, and using induction,
there exists a neighborhood $I$ of $0$ in $\R$ and an n-plaque
$p\colon I^k\to X$ such that $[p(te_i)] = [p_i(t)]$, therefore
$$\eqalign{ \omega_{p(0})(e_1, \cdots e_k) &= 
\omega_{p(0)}([p(te_1)], \cdots ,[p(te_k)]) \cr
&= \omega_{p(0)}([p_1(t)], \cdots ,[p_k(t)]) \cr
&= \omega_F(v_1, \cdots ,v_k) \neq 0 \cr}
$$
Therefore $\omega$ is 1-1.
Now we prove that $\Psi$ is onto: Let $(\omega_p)_{p\in P(X)} \in 
\epsilon_3^k(X)$, define $\underline{\omega}$ as follows: let $F\in X$
and let $v_1, \cdots ,v_k\in \underline{T_FX}$, let $p\in P(X)$ such
that $[p(te_i)]= v_i$ (exists because $F$ is not strongly transverse),
define
$$\underline{\omega}_F(v_1, \cdots ,v_k) = (\omega_p)(0)(e_1, \cdots ,e_k)$$
This value is independent of $p$ because if $p^{\prime}$ is another plaque
as before then
$$(\omega_p)(0)(e_1, \cdots ,e_k)= (\omega_{p^{\prime}}(0)(e_1, \cdots e_k)$$
because of the tangent condition. Next we prove that $\underline{\omega}$
is smooth: Let $\xi_1, \cdots ,\xi_k \in \Gamma_I(X)$ and let
$p\colon U\subset \R^n \to X$ be a plaque. We shall prove that
$$\underline{\omega}(p(r)(\xi_1(p(r), \cdots \xi_k(p(r))\colon U \to \R$$
is smooth. Let $q_1, \cdots , q_n$ be $(n+1)$-plaques such that
$[q_i(r,t)]_t = \xi_i(p(r)) \forall r$. Let
$q\colon V\subset \R^{n+k} \to X$ be a plaque as in the definition of 
transverse point such that
$$p(r,0, \cdots ,t_i, 0, \cdots ,0)= q_i(r,t_i$$
then $q(r, 0, \cdots ,0)= p(r), \forall r$ and 
$$[q(r,0, \cdots ,0,t_i, 0, \cdots ,0)]_{t_i} = [q_i(r,t_i)]_{t_i}$$
$$\eqalign{& \underline{\omega}_{p(r)}(\xi_1(p(r)), \cdots ,\xi_k(p(r))\cr
&= \underline{\omega}_{q(r,0)}([q_1(r,t)]_t, \cdots ,
[q_k(r,t)]_t) \cr
&= (\omega_q)(r,0)([(r,t_1,0, \cdots ,0)]_{t_1}, 
[(r,0,t_2, \cdots ,0)]_{t_2}, \cdots ,[(r,\cdots ,t_k)]_{t_k} \cr
&= (\omega_q)(r,0)(v_1(r), \cdots ,v_k(r)) \cr}$$
and $v_i(r)= [(r,0, \cdots ,0, t_i, 0, \cdots ,0)]_{t_i}$ is a smooth
vector field on $U\subset \bar{U}$, therefore $\underline{\omega}$ is
smooth. Notice that if $p\colon U\subset \R^n \to X$ is a plaque then
$$\eqalign{(p^*\underline{\omega)}(v_1, \cdots ,v_k) &=
\underline{\omega}_{p(r)}([p(tv_1)], \cdots ,[p(tv_k)]) \cr
&= (\omega_p)(r)(v_1, \cdots ,v_k) \cr}$$
It follows from the definition that $\Psi$ is onto since
$\Psi (\underline{\omega})= \omega$ and this proves the proposition.

\declare Remarks 3.1.

\item{1.}
The theorem remains still remains valid under a weaker condition than
non transversality: for every pair of $(n+1)$-plaques
$p_1, p_2\colon U\times V \subset \R^n \times \R \to X$ such that
$p_1(r,0)= p_2(r,0), \forall r\in U$ there exists a $(n+2)$-plaque
$$q\colon U^{\prime}\times V^{\prime}\times \times V^{\prime \prime}
\subset \R^n \to \R \times \R \to X$$ 
such that
$$[q(r,t_1,0)]_{t_1} = [p_1(r,t_1)]_{t_1} \ \ \forall r\in U^{\prime}$$
$$[q(r,0,t_2)]_{t_2} = [p_2(r,t_2)]_{t_2} \ \ \forall r\in U^{\prime}$$

\item{2.}
For any $X$ we have $\epsilon_1^k(X) \subset \epsilon_2^k(X)$.

\item{3.}
In general $\epsilon_2^k(X)$ is larger than $\epsilon_1^k(X)$:
consider the example
$$X= \{ (x,0) / x\in \R \} \cup \{ (0,y) / y\in \R \}$$
where $P(X)$ is the set of smooth maps $p\colon U\subset \R^n \to X$
whose image is a subset of the x-axis only, or a subset of the y-axis only.
Then $c^{\infty}(X)$ is the set of pairs $(f_1, f_2)$ such that 
$f_1, f_2\colon \R \to \R$ is smooth and $f_1(0)= f_2(0)$. With the standard
diffeology $T_{(0,0)}X = \R \sqcup \R$ and $\underline{T_{(0,0)}}X= \{ 0 \}$,
otherwise, $T_FX = \R= \underline{T_F}(X)$. The point $(0,0)$ is transverse
and $\Gamma(X)$ is the set of pairs 
$(f_1{\partial \over \partial x}, f_2{\partial \over \partial y})$ such
that $f_1, f_2\colon \R \to \R$ is smooth and $f_1(0)= f_2(0) = 0$.
Let 
$$\xi_1(x,y)= x{\partial \over \partial x} \ \ \ \xi_2(x,y) = 
y{\partial \over \partial y}$$
Since every smooth function $h\colon \R \to \R$ such that $h(0)= 0$
satisfies $h(t) = tg(t)$ where $g$ is also smooth then every field
is expressed as 
$$\xi(x,,y)= h_1(x)\xi_1(x,y) + h_2(y)\xi_2(x,y)$$
where $h_1, h_2\colon \R \to \R$ and conversely, define
$$\omega(\xi(x,y)) = h_1(x)(x^2+1) + h_2(y)(2y^2+1)$$
then $\omega \in \epsilon_2(X)$ and satisfies
$\omega(\xi_1(0,0)= 1$, however $\xi_1(0,0)= 0$, therefore there is not
$\underline{\omega}\in \epsilon_1(X)$ which correspond to $\omega$.
Similar constructions may be given with any $k>1$.

\item{4.}
A construction similar to the above may be done for any $X$ such that
$\Gamma(X)$ is a free module and there exists a base $\{ \xi_i / i\in I\}$
such that at some point $F$ the set of vectors 
$$\{ \xi_i(F) / i\in I \} \subset T_FX$$
is linearly dependent (this occurs if $dim T_FX \neq Card I$ for some $F$).
This occurs if $\Gamma(X)$ is free and $X$ has transverse points.

\item{5.}
$\epsilon_1(X) = \epsilon_2(X)$ if $\Gamma(X)$ is a free module and 
there exists a base $\{ \xi_i / i\in I\}$ of $\Gamma(X)$ such that
$\{ \xi_i(F) / i\in I \}$
is linearly independent $\forall F\in X$ (this occurs for example with
$\R^M, M\in \N \cup \{ \infty \}$ and manifolds with trivial tangent
bundles such as $S^1$), because if
$\omega(\xi_i) = h_i$ then
$$\omega(\xi = \sum_if_i\xi_i) = \sum_i h_if_i$$
and $\xi(F)=0$ implies $f_i(F)= 0, \forall i\in I$

\vfill
\eject

\varsection  References


\refno \DonatoGeometrie.
P. Donato,
G\'eom\'etrie des orbites coadjointes des groupes de
dif\-feo\-mor\-phismes,
in: C. Albert, ed., {\it G\'eom\'etrie Symplectique et M\'ecanique},
Lecture Notes in Mathematics {\bf 1416} (Springer, Berlin, 1988)
84--104.

\refno \DonatoVirasoro.
P. Donato,
Les diffeomorphismes du circle comme orbit symplectique dans les
moments de Virasoro,
Preprint CPT--92/P.2681, CNRS--Luminy, 1992.

\refno \SouriauAlgoritme.
J. M. Souriau,
Un algoritme g\'en\'erateur de structures quantiques,
Ast\'erisque, hors s\'erie (1985) 341--399.

\refno \Sou.
J. M. Souriau,
Groupes diff'rentiels,
Lecture notes in math., 838 (1980) 91-128.

\refno \TorreTangen.
C. Torre,
A Tangent bundle on diffeological spaces,
Math/9801046 (1998).

\refno \Torre.
C. Torre,
Examples of smooth diffeological spaces,
Preprint (1999).

\refno \TorreBan.
C. Torre and A. Banyaga,
A symplectic structure on coadjoint orbits of diffeomorphism subgroups,
Ciencia y Tecnolog¡a 17 No 2, (1993) 1--14.

\bye